\documentclass[11pt, twoside]{amsart}
\usepackage{amsmath, amsthm, amscd, amsfonts, amssymb, graphicx, color, delarray}
\usepackage[bookmarksnumbered, plainpages, backref]{hyperref}

\begin{document}

\title{Integrated and Differentiated Spaces of Triangular Fuzzy Numbers}

\author{Murat Kiri\c{s}ci}

\address{[Murat Kiri\c{s}ci] Department of Mathematical Education, Hasan Ali Y\"{u}cel Education Faculty,
Istanbul University, Vefa, 34470, Fatih, Istanbul, Turkey \vskip 0.1cm }
\email{mkirisci@hotmail.com, murat.kirisci@istanbul.edu.tr}

\begin{abstract}
Fuzzy sets are the cornerstone of a non-additive uncertainty theory, namely possibility theory,
 and of a versatile tool for both linguistic and numerical modeling. Numerous works now combine
 fuzzy concepts with other scientific disciplines as well as modern technologies. In mathematics
 fuzzy sets have triggered new research topics in connection
with category theory, topology, algebra, analysis. In this paper,
 we use the triangular fuzzy numbers for matrix domains of sequence spaces with infinite matrices.
 We construct the new space with triangular fuzzy numbers and investigate to structural, topological and algebraic
 properties of these spaces.
\end{abstract}

\maketitle

\pagestyle{plain} \makeatletter
\theoremstyle{plain}
\newtheorem{thm}{Theorem}[section]
\numberwithin{equation}{section}
\numberwithin{figure}{section}  
\theoremstyle{plain}
\newtheorem{pr}[thm]{Proposition}
\theoremstyle{plain}
\newtheorem{exmp}[thm]{Example}
\theoremstyle{plain}
\newtheorem{cor}[thm]{Corollary} 
\theoremstyle{plain}
\newtheorem{defin}[thm]{Definition}
\theoremstyle{plain}
\newtheorem{lem}[thm]{Lemma} 
\theoremstyle{plain}
\newtheorem{rem}[thm]{Remark}
\numberwithin{equation}{section}

\section{Introduction}

The concept of fuzziness was first discovered and introduced by Zadeh \cite{Zadeh}. Between 1965 and 1975, Zadeh broadened the
foundation of the fuzzy set theory by establishing fuzzy
similarity relations, linguistic hedges, and fuzzy
decision-making.\\

To understand fuzzy sets it is useful to first consider some of the fundamentals of classical set theory. In particular, consider the two fundamental laws of Boolean algebra the law of excluded middle and law of contradiction. In logic, the proposition every proposition is either true or false excludes any third, or middle, possibility, which gave this principle the name of  the law of excluded middle. From this we derive that a prediction may be only true or false. Thus, in classical, Boolean or crisp, set theory, membership of an element  $x$  in a set  $A$ , is defined by a characteristic function which assigns a value of either $1$ or $0$ to each individual in the universal set  $X$. This is how it discriminates between members and nonmembers of the crisp set  $A$.\\

The applications of fuzziness are uncountable and varied. The fuzzy
literature is vast and huge and goes to all possible directions. Fuzzy sets were taken up with enthusiasm by engineers, computer scientists and operations researchers. A notable reason is the relationship that fuzzy sets have with multivalued logic, offering decision possibilities such as "maybe true" and "maybe false", suitably quantified, in addition to the traditional dichotomy of "true" and "false". Vast, interconnected arrays of such decisions can be processed rapidly by computers.\\

In classical, or crisp, sets the transition for an element in the universe between membership and nonmembership in a given set is abrupt and well defined.
For an element in a universe that contains fuzzy sets, this transition can be gradual. This transition among various degrees of membership can be thought of as conforming to the fact that the boundaries of the fuzzy sets are vague and ambiguous. Hence, membership of an element from the universe in this set is measured by a function that attempts to describe vagueness and ambiguity.\\

The classes of bounded and convergent sequences of fuzzy numbers with respect to the Hausdorff metric are given by Matloka\cite{matloka}.
The spaces of bounded and convergent sequences of fuzzy numbers was studied also by Nanda\cite{nanda} and was shown that these spaces are complete metric spaces.
Talo and Basar\cite{TaBa0} have defined the space $bv_{p}(F)$ of fuzzy numbers which is an extension of the space $bv_{p}$. In \cite{TaBa1}, dual spaces of the classical sets of sequences of fuzzy numbers was computed and characterized some matrix classes. The quasilinearity of the classical sets of sequences of fuzzy numbers have been studied in \cite{TaBa2}. Quite recently, the matrix domains of the sequence spaces with fuzzy numbers are studied\cite{seng1, seng2, TaCa, zararsiz0, zararsiz1, zararsizSeng}.\\

In this paper, we construct the integrated and differentiated spaces using the triangular fuzzy numbers and compute the dual spaces of new spaces.
Finally, we characterize some matrix classes.

\section{Preliminaries, Background and Notation}

\subsection{Definitions and Notions}
First, we give $\Omega$ and $\Gamma$ matrices.\\

 The matrix $\Omega=(a_{nk})$ defined by $a_{nk}=k$, $(1\leq k \leq n)$ and $a_{nk}=0$, $(k>n)$, and the matrix $\Gamma=(b_{nk})$ defined by
$b_{nk}=1/k$, $(1\leq k \leq n)$ and $b_{nk}=0$, $(k>n)$, respectively, i.e.,

\begin{eqnarray*}
a_{nk}= \left[ \begin{array}{cccccc}
1 & 0 & 0 & 0 & \cdots \\
1 & 2 & 0 & 0 & \cdots  \\
1 & 2 & 3 & 0 & \cdots \\
1 & 2 & 3 & 4 & \cdots \\
\vdots & \vdots &  \vdots&  \vdots& \ddots
\end{array} \right]
\quad \textrm{ and } \quad
b_{nk}= \left[ \begin{array}{cccccc}
1 & 0 & 0 & 0 & \cdots \\
1 & 1/2 & 0 & 0 & \cdots  \\
1 & 1/2 & 1/3 & 0 & \cdots \\
1 & 1/2 & 1/3 & 1/4 & \cdots \\
\vdots & \vdots &  \vdots&  \vdots& \ddots
\end{array} \right]
\end{eqnarray*}

We can give the matrices $\Omega^{-1}=(c_{nk})$ and $\Gamma^{-1}=(d_{nk})$ which are inverse of the above matrices by $c_{nk}=1/n$, $(n=k)$, $c_{nk}=-1/n$, $(n-1=k)$, $c_{nk}=0$, $(other)$ and $d_{nk}=n$, $(n=k)$, $d_{nk}=-n$, $(n-1=k)$, $d_{nk}=0$, $(other)$, respectively, i.e.,

\begin{eqnarray*}
c_{nk}= \left[ \begin{array}{cccccc}
1 & 0 & 0 & 0 & \cdots \\
-1/2 & 1/2 & 0 & 0 & \cdots  \\
0 & -1/3 & 1/3 & 0 & \cdots \\
0 & 0 & -1/4 & 1/4 & \cdots \\
\vdots & \vdots &  \vdots&  \vdots& \ddots
\end{array} \right]
\quad \textrm{ and } \quad
d_{nk}= \left[ \begin{array}{cccccc}
1 & 0 & 0 & 0 & \cdots \\
-2 & 2 & 0 & 0 & \cdots  \\
0 & -3 & 3 & 0 & \cdots \\
0 & 0 & -4 & 4 & \cdots \\
\vdots & \vdots &  \vdots&  \vdots& \ddots
\end{array} \right]
\end{eqnarray*}

The set of all sequences denotes with $\omega := \mathbb{C}^{\mathbb{N}}:=\{x=(x_{k}): x: \mathbb{N}\rightarrow \mathbb{C}, k\rightarrow x_{k}:=x(k)\}$, where $\mathbb{C}$ denotes the complex field and $\mathbb{N}=\{0,1,2,\ldots\}$. Each linear subspace of $\omega$ (with the induced addition and scalar multiplication) is called a \emph{sequence space}. We will write $\phi, \ell_{\infty}, c ~\textrm{ and }~ c_{0}$ for the sets of all finite, bounded, convergent and null sequences, respectively. It obviously that these sets are subsets of $\omega$.\\

The concept of integrated and differentiated spaces was employed as
\begin{eqnarray*}
\int X=\left\{x=(x_{k})\in\omega: (kx_{k})\in X\right\} \quad ~ \textrm{and}~ \quad d(X)=\left\{x=(x_{k})\in\omega: (k^{-1}x_{k})\in X\right\},
\end{eqnarray*}
where $X$ is a sequence space. Integrated and differentiated sequence spaces which emerged from the concept of the sequence space $\ell_{1}$. The integrated and differentiated
 sequence spaces were initiated by Goes and Goes \cite{Goes}.  Subramanian et al. \cite{SubRaoG} gave the integrated rate space and studied some properties. Until Kiri\c{s}ci paper\cite{kirisci1}, there was no work on differentiated sequence spaces.  Kiri\c{s}ci\cite{kirisci1} defined new integrated and differentiated spaces by using the new $\Omega$ and  $\Gamma$ matrices and investigated some properties. In\cite{kirisci2}, integrated and differentiated sequence spaces are defined with Riesz mean.\\

Now, we give some basic definitions and notions related to the fuzzy numbers.\\

The set of all bounded and closed intervals on $\mathbb{R}$ are defined by $E_{i}=\{a=[a^{-},a^{+}]:~a^{-}\leq x \leq a^{+}, ~a^{-} ~ \textrm{and}~ a^{+}\in \mathbb{R}\}$. If we define the metric
\begin{eqnarray}\label{metric}
d(a,b)=\max\left\{\left|a^{-}-b^{-}\right|, \left|a^{+}-b^{+}\right|\right\},
\end{eqnarray}
for $a,b\in E_{i}$, then, the pair $(E_{i},d)$ is a complete metric space. The fuzzy subset of $X$ which is nonempty set, is a nonempty subset $\{(x,S(x)): ~ x\in X\}$ of $X\times [0,1]$, where $S:X\rightarrow [0,1]$. If we take $X=\mathbb{R}$, then the function $S:\mathbb{R}\rightarrow[0,1]$ is a subset  of nonempty base space $\mathbb{R}$.\\

We give the following properties:\\
\textbf{i.} There exists an $x_{0}\in \mathbb{R}$ such that $S(x_{0})=1$.\\
\textbf{ii.} $S[\lambda x+(1-\lambda)y]\geq \min\{S(x),S(y)\}$, for any $x,y\in \mathbb{R}$ and $\lambda\in[0,1]$.\\
\textbf{iii.} $S$ is upper semicontinuous.\\
\textbf{iv.} The closure of $\{x\in\mathbb{R}: ~S(x)>0 \}$ is compact(denoted by $S^{0}$).\\

If the function $S$ provides the conditions (i)-(iv), then $S$ is called the fuzzy number. The $\alpha-$level set of fuzzy set $S$
on $\mathbb{R}$ is defined as $S(\alpha)=\{x\in \mathbb{R}: ~ S(x)\geq \alpha\}$ for each $\alpha\in (0,1]$.\\

If we choose the membership function of the triangular fuzzy number $u_{(t_{1},t_{2})}$,
\begin{eqnarray}\label{triangular}
u_{(t_{1},t_{2})}=\left\{\begin{array}{ccl}
\big(x-(u-t_{1})\big)t_{1}^{-1}&, & x\in [u-t_{1},u]\\
\big((u+t_{2})-x\big)t_{2}^{-1}&, & x\in [u,u+t_{2}]\\
0&, & otherwise
\end{array}\right.
\end{eqnarray}
then, for $t_{1}\leq t_{2}$, $(t_{1},t_{2}\in\mathbb{R})$, $u_{(t_{1},t_{2})}$ can be represented with the notation $u_{(t_{1},t_{2})}=(u-t_{1},u,u+t_{2})$.
$u_{(t_{1},t_{2})}$ is a real number, when $t_{1}=t_{2}=0$. For convenience, we will understand the $(t_{1},t_{2})-$type fuzzy numbers, when we say triangular fuzzy numbers, through the text. We denote the set $F$ as
\begin{eqnarray}\label{triangularset}
F=\big\{(u-t_{1},u,u+t_{2}): ~ t_{1},t_{2}\in\mathbb{R}, ~ t_{1}\leq t_{2}\big\}.
\end{eqnarray}

The points $u-t_{1},u,u+t_{2}$  of the triangular fuzzy number $u_{(t_{1},t_{2})}$ are called first, middle, end points, respectively. The notation $u$
means that the height of the fuzzy number $u_{(t_{1},t_{2})}$ is $1$ at the point $u$. There is no unique set of fuzzy numbers in the set in form (\ref{triangularset}). Conversely, there are infinitely-many sets of fuzzy numbers which are different from each other according to structure of their elements.\\

\subsection{$\oplus$, $\ominus$, $\otimes$, $\oslash$ Operators}

Errors occur according to algebraic operations, when we choose the fuzzy numbers with $\alpha-$cut sets.
We can explain this with an example. Take any fuzzy number $a=[a^{-}(\alpha),a^{+}(\alpha)]$. For $\alpha\in[0,1]$,
$a-a=[a^{-}(\alpha),a^{+}(\alpha)]-[a^{-}(\alpha),a^{+}(\alpha)]=[a^{-}(\alpha)-a^{+}(\alpha),a^{+}(\alpha)-a^{-}(\alpha)]$.
It is clear that this equation is nonzero as expected in the classical mean. Therefore, $\theta$ denotes the fuzzy zero.

We define the algebraic structure of the set of all $(t_{1},t_{2})-$ type fuzzy number.
Take $u_{(t_{1},t_{2})}\in F$ and choose $\lambda \in \mathbb{R}$. If $v_{(t_{1},t_{2})}$ is nonzero fuzzy number, then,
\begin{eqnarray}\label{oplus}
u_{(t_{1},t_{2})}\oplus v_{(t_{1},t_{2})}=(u-t_{1},u,u+t_{2})\oplus (v-t_{1},v,v+t_{2})
\end{eqnarray}
\begin{eqnarray*}
=(u+v-t_{1},u+v,u+v+t_{2})
\end{eqnarray*}

\begin{eqnarray}\label{ominus}
u_{(t_{1},t_{2})}\ominus v_{(t_{1},t_{2})}=(u-t_{1},u,u+t_{2})\ominus (v-t_{1},v,v+t_{2})
\end{eqnarray}
\begin{eqnarray*}
=(u-v-t_{1},u-v,u-v+t_{2})
\end{eqnarray*}

\begin{eqnarray}\label{otimes}
u_{(t_{1},t_{2})}\otimes v_{(t_{1},t_{2})}=(u-t_{1},u,u+t_{2})\otimes (v-t_{1},v,v+t_{2})
\end{eqnarray}
\begin{eqnarray*}
=(uv-t_{1},uv,uv+t_{2})
\end{eqnarray*}

\begin{eqnarray}\label{oslash}
u_{(t_{1},t_{2})}\oslash v_{(t_{1},t_{2})}=(u-t_{1},u,u+t_{2})\oslash (v-t_{1},v,v+t_{2})
\end{eqnarray}
\begin{eqnarray*}
=(u\div v-t_{1},u\div v,u\div v+t_{2})
\end{eqnarray*}

\begin{eqnarray}\label{oscalar}
\alpha \otimes u_{(t_{1},t_{2})}=(\alpha u-t_{1},\alpha u,\alpha u+t_{2})
\end{eqnarray}
\begin{eqnarray*}
=(r-t_{1},r,r+t_{2})=r_{(t_{1},t_{2})}.
\end{eqnarray*}

We find the identity element of $F$. We choose $u_{(t_{1},t_{2})}=(u-t_{1},u,u+t_{2})$ in $F$.
Then, we obtain
$(u-t_{1},u,u+t_{2})+(0-t_{1},0,0+t_{2})=(u+0-t_{1},u+0,u+0+t_{2})$. It means that $\theta_{(t_{1},t_{2})}=(0-t_{1},0,0+t_{2})$
is considered as the identity element of $F$ according to operation which is given in (\ref{oplus}).\\

Now, we give the inverse of triangular fuzzy number according to addition. The inverse of fuzzy number $u_{(t_{1},t_{2})}\in F$
is equal to $-u_{(t_{1},t_{2})}=(-u-t_{1},-u,-u+t_{2})$ and $-u_{(t_{1},t_{2})}$ is a fuzzy number. Therefore, we can write
$x_{(t_{1},t_{2})}+u_{(t_{1},t_{2})}=\theta_{(t_{1},t_{2})}$, where $\theta=\theta_{(t_{1},t_{2})}=(0-t_{1},0,0+t_{2})$ denotes the fuzzy zeros of the sets $F$.
From this idea, we can see that the fuzzy zero of each element of the set $F$ is different. It is clear that the representation $u_{(t_{1},t_{2})}=(u-t_{1},u,u+t_{2})$ is unique. Then, for every $-u_{(t_{1},t_{2})}=(-u-t_{1},-u,-u+t_{2})$,
certainly $u_{(t_{1},t_{2})}$ is unique.

\begin{lem}
Let $t_{1}\leq t_{2}$ $(t_{1},t_{2}\in\mathbb{R})$. All sets of the form $F$ are linear spaces according to algebraic operations (\ref{oplus}) and (\ref{oscalar}).
\end{lem}

\subsection{Topological structure of the set $F$}

Define the metric $\overline{d}:F\times F\rightarrow \mathbb{R}$ as follows:
\begin{eqnarray*}
\overline{d}\left(u_{(t_{1},t_{2})},v_{(t_{1},t_{2})}\right):=\max\big\{|u-v-t_{1}|, |u-v|, |u-v+t_{2}|\big\}.
\end{eqnarray*}
It is clear that the $(F, \overline{d})$ is a complete metric space.\\

Applications have shown that the spread of fuzziness should not be very large. Then, it is necessary to make the value $\max\{|u-v-t_{1}|, |u-v|, |u-v+t_{2}|\}$
as small as possible. In fact, we theoretically know that this is not necessary, but this is absolutely necessary in practice. Let us explain this thought with an example: Define the expression "approximately 5" as $5_{(t_{1},t_{2})}=(-4-t_{1},5,15+t_{2})$. In applications, this expression is taken as $5_{(t_{1},t_{2})}=(5-t_{1},5,5+t_{2})$, $(0\leq t_{1}\leq t_{2}<1)$ and this choice is more accurate than $5_{(t_{1},t_{2})}=(-4-t_{1},5,15+t_{2})$.\\

It should be noted that a fuzzy number is determined according to specific processes. For example, we take specific systems $X$ and $Y$ which are different.
For system $X$, choose $S$ as "approximately 5". Then, we comprehend that the expression "approximately 5" may not in the same sense for another system $Y$.
Naturally, algebraic properties of these systems are different. We also give an example to this situation: Assume that the spread of left and right fuzziness of every number $S$ be equal to $T$ in the system $X$. Then, fuzzy zero is equal to $(-2t, 0, 2t)$, $(0\leq t<1)$ for the system $X$ and this fuzzy zero is unique for $X$.

\begin{thm}
The matrices $\Omega$ and $\Gamma$ are regular.
\end{thm}

\begin{proof}
Let $\widehat{u}=(u_{(t_{1},t_{2})}^{k})$ be a sequence of triangular fuzzy numbers. We must show that
if for $n\rightarrow \infty$, $\bar{d}\left(u_{(t_{1},t_{2})}^{k}, u_{(t_{1},t_{2})}^{0}\right)\rightarrow 0$,
then for $n\rightarrow \infty$, $\bar{d}\left(\sum_{k=1}^{n}\left|ku_{(t_{1},t_{2})}^{k}\right|, u_{(t_{1},t_{2})}^{0}\right)\rightarrow 0$.
Suppose that for $n\rightarrow \infty$, $\bar{d}\left(u_{(t_{1},t_{2})}^{k}, u_{(t_{1},t_{2})}^{0}\right)\rightarrow 0$ and choose
$\varepsilon >0$. Then, there exist a positive integer $N$ such that $\bar{d}\left(u_{(t_{1},t_{2})}^{k}, u_{(t_{1},t_{2})}^{0}\right)<\varepsilon$
for $n\geq N$. Then, for $n\geq N$ and $N \in \mathbb{N}$, $\bar{d}\left(\Omega u_{(t_{1},t_{2})}^{k}, u_{(t_{1},t_{2})}^{0}\right)= \bar{d}\left(\sum_{k=1}^{n}\left|ku_{(t_{1},t_{2})}^{k}\right|, u_{(t_{1},t_{2})}^{0}\right)<\varepsilon$. Therefore the matrix $\Omega$ is regular.\\

Similarly, we can show that the matrix $\Gamma$ is regular.
\end{proof}

\subsection{Classical sequence spaces of fuzzy number}

Define the $f:\mathbb{N}\rightarrow F$, $k\rightarrow f(k)=u_{(t_{1},t_{2})}^{k}$.
The function $f$ is defined as a sequence of triangular fuzzy numbers and is represented $\widehat{u}=\big(u_{(t_{1},t_{2})}^{k}\big)$.
The set of all sequences of triangular fuzzy numbers denotes with $\omega(F)=\bigg\{\widehat{u}=\big(u_{(t_{1},t_{2})}^{k}\big): ~ u:\mathbb{N}\rightarrow F, ~ u(k)=\big(u_{(t_{1},t_{2})}^{k}\big)=\big(u^{k}-t_{1},u^{k},u^{k}+t_{2}\big)\bigg\}$, where $u^{k}-t_{1}\leq u^{k} \leq u^{k}+t_{2}$, $(t_{1},t_{2}\in \mathbb{R})$ and $u_{(t_{1},t_{2})}^{k}\in F$ for all $k\in \mathbb{N}$. Each subspace of $\omega(F)$ is called a sequence space of fuzzy numbers. In this place, the elements $u^{k}-t_{1}$, $u^{k}$, $u^{k}+t_{2}$ is expressed first, middle, end points of general term of a sequences of fuzzy numbers, respectively. If degree of membership at $u^{k}$ is equal to $1$, then $\widehat{u}$ is a $(t_{1},t_{2})-$type fuzzy number, if it is not equal to $1$, then $\big(u_{(t_{1},t_{2})}^{k}\big)$ is a sequence of the fuzzy sets.\\

We define the bounded, convergent, null and absolutely $p-$summable sequence spaces as follows\cite{zararsiz0, zararsiz1}:\\
\begin{eqnarray*}
&&\ell_{\infty}(F)=\left\{\widehat{u}=\big(u_{(t_{1},t_{2})}^{k}\big)\in\omega(F): \sup_{k\in\mathbb{N}}\overline{d}\left(u_{(t_{1},t_{2})}^{k}, \theta\right)<\infty\right\},\\
&&c(F)=\left\{\widehat{u}=\big(u_{(t_{1},t_{2})}^{k}\big)\in\omega(F): \lim_{k}\overline{d}\left(u_{(t_{1},t_{2})}^{k}, u_{(t_{1},t_{2})}^{0}\right)=0, \quad u_{(t_{1},t_{2})}^{0}\in F\right\},\\
&&c_{0}(F)=\left\{\widehat{u}=\big(u_{(t_{1},t_{2})}^{k}\big)\in\omega(F): \lim_{k}\overline{d}\left(u_{(t_{1},t_{2})}^{k}, \theta\right)=0\right\},\\
\end{eqnarray*}
\begin{eqnarray*}
&&\ell_{p}(F)=\left\{\widehat{u}=\big(u_{(t_{1},t_{2})}^{k}\big)\in\omega(F): \sum_{k}\overline{d}\left(u_{(t_{1},t_{2})}^{k}, \theta\right)^{p}<\infty, ~1\leq p<\infty \right\}.
\end{eqnarray*}

Let $X(F)\subset \omega(F)$ and define the function $\|.\|:X(F)\rightarrow \mathbb{R}$. Suppose that the function $\|.\|$ is satisfied the following conditions:\\
i. $\|\widehat{u}\|=0 \Leftrightarrow \widehat{u}=\theta$,\\
ii. $\|\alpha\widehat{u}\|=|\alpha|\|\widehat{u}\|, \quad \alpha\in \mathbb{R}$,\\
iii. $\|\widehat{u}+\widehat{v}\|\leq \|\widehat{u}\|+\|\widehat{v}\|$.\\

Then, the function $\|.\|$ is called norm and $X(F)$ is called normed sequence space of the $\big(t_{1},t_{2}\big)-$fuzzy numbers. If $X(F)$ is complete with respect to the norm $\|.\|$, then $X(F)$ is called complete normed sequence space of the $\big(t_{1},t_{2}\big)-$fuzzy numbers.\\

\begin{lem}(\cite{zararsiz0})
The bounded, convergent, null sequence spaces are complete normed sequence spaces with the norm defined by
\begin{eqnarray*}
\|\widehat{u}\|=\sup_{k}\max\left\{|u^{k}-v^{k}-t_{1}|, |u^{k}-v^{k}|, |u^{k}-v^{k}+t_{2}|\right\},
\end{eqnarray*}
where $\widehat{u}$ is in the any sets of $\{\ell_{\infty}(F), c(F), c_{0}(F)\}$.
\end{lem}

Let $A=(a_{nk})$ be an infinite matrix and choose two spaces of triangular fuzzy valued sequences $X(F)$ and $Y(F)$. Then, we define
the \emph{real} matrix mapping from $X(F)$ to $Y(F)$ as $A: X(F) \rightarrow Y(F)$. That is, if we take $\widehat{u}=\big(u_{(t_{1},t_{2})}^{k}\big)\in X(F)$,
then we can write $A\widehat{u}=\big[(Au_{(t_{1},t_{2})})^{n}\big]\in Y(F)$, where,
\begin{eqnarray}\label{matrixF}
&&\left(Au_{(t_{1},t_{2})}\right)^{n}=\sum_{k}a_{nk}u_{(t_{1},t_{2})}^{k}=\left(\sum_{k}a_{nk}u^{k}-t_{1}, \sum_{k}a_{nk}u^{k}, \sum_{k}a_{nk}u^{k}+t_{2}\right).
\end{eqnarray}
In (\ref{matrixF}), the series $\sum_{k}a_{nk}u^{k}-t_{1}$, $\sum_{k}a_{nk}u^{k}$, $\sum_{k}a_{nk}u^{k}+t_{2}$ are convergent for all $n\in \mathbb{N}$.
That is, $A: X(F) \rightarrow Y(F)$ if and only if for $u\in X(F)$ and all $n\in \mathbb{N}$, the series of $\sum_{k}a_{nk}u^{k}-t_{1}$, $\sum_{k}a_{nk}u^{k}$, $\sum_{k}a_{nk}u^{k}+t_{2}$ are convergent.\\

The domain of an infinite matrix $A$ in $X(F)$ is defined by

\begin{eqnarray*}
[X(F)]_{A}=\left\{(u_{(t_{1},t_{2})}^{k})\in \omega(F): Au_{(t_{1},t_{2})}^{n}\in X(F) \right\}.
\end{eqnarray*}

\section{New spaces with triangular fuzzy numbers}
Let $u_{t_{1}}^{k}\leq u^{k} \leq u_{t_{2}}^{k}$ and $u=(u_{(t_{1},t_{2})}^{k}\in F$ for all $k\in \mathbb{N}$.\\

Now, we introduce the integrated spaces with triangular fuzzy numbers by
\begin{eqnarray*}
[X(F)]_{\Omega}=\left\{\widehat{u}=\left(u_{(t_{1},t_{2})}^{k}\right)\in \omega(F): \Omega u_{(t_{1},t_{2})}^{k} \in X(F)\right\}
\end{eqnarray*}
and the differentiated spaces with triangular fuzzy numbers by
\begin{eqnarray*}
[X(F)]_{\Gamma}=\left\{\widehat{u}=\left(u_{(t_{1},t_{2})}^{k}\right)\in \omega(F): \Gamma u_{(t_{1},t_{2})}^{k} \in X(F)\right\}
\end{eqnarray*}
where $X=\{\ell_{\infty}, c, c_{0}\}$.\\

Let us define the sequence of fuzzy numbers $v=\left(v_{(t_{1},t_{2})}^{k}\right)$ and $y=\left(y_{(t_{1},t_{2})}^{k}\right)$, as the $\Omega-$transform and $\Gamma-$tarnsform of a sequence
of fuzzy numbers $u=\left(u_{(t_{1},t_{2})}^{k}\right)$, respectively; that is, for $k,n\in \mathbb{N}$,
\begin{eqnarray}\label{transfseq1}
&&v_{(t_{1},t_{2})}^{n}=\sum_{k=1}^{n}\left|ku_{(t_{1},t_{2})}^{k}\right|
\end{eqnarray}
\begin{eqnarray}\label{transfseq2}
&&y_{(t_{1},t_{2})}^{n}=\sum_{k=1}^{n}\left|\frac{1}{k}u_{(t_{1},t_{2})}^{k}\right|
\end{eqnarray}
where $u_{(t_{1},t_{2})}^{-1}=\theta$.

\begin{thm}\label{isotheo}
Integrated and differentiated spaces derived by triangular fuzzy numbers are norm isomorphic to the
classical sets consisting of the bounded, convergent and null sequences of fuzzy numbers.
\end{thm}

\begin{proof}
We will show that there is a linear isometry between integrated bounded space with triangular fuzzy numbers and bounded sequence space of fuzzy numbers.
We consider the transformation defined $\Phi$, from $[\ell_{\infty}(F)]_{\Omega}$ to $\ell_{\infty}(F)$ by $u \mapsto v =\Phi u = \sum_{k=1}^{n}\left|ku_{(t_{1},t_{2})}^{k}\right|$. Then, it is
clear that the equality $\Phi(u+v)=\Phi(u)+\Phi(v)$ is holds. Choose $\lambda \in \mathbb{R}$. Then,
\begin{eqnarray*}
\Phi(\lambda u)=\Phi(\lambda u_{(t_{1},t_{2})}^{k})&=& \sum_{k=1}^{n}\left|\lambda ku_{(t_{1},t_{2})}^{k}\right|\\
&=& \lambda \sum_{k=1}^{n}\left|ku_{(t_{1},t_{2})}^{k}\right|=\lambda \Phi u.
\end{eqnarray*}
Therefore, we can say that $\Phi$ is linear.\\

We take $v\in \ell_{\infty}(F)$ and define the sequence $u$ such that $u=(u_{(t_{1},t_{2})}^{k})=(c_{nk}v_{(t_{1},t_{2})}^{k})$, where $(c_{nk})$  is an inverse of the matrix $\Omega$.
\begin{eqnarray*}
\|u\|_{[\ell_{\infty}(F)]_{\Omega}}&=& \sup_{k} \bar{d} \left(\Omega u_{(t_{1},t_{2})}^{k}, \theta \right)\\
&=&\sup_{k} \bar{d}\left(v_{(t_{1},t_{2})}^{k}, \theta \right)=\|v\|_{\ell_{\infty}(F)}
\end{eqnarray*}
Therefore, we can say that $\Phi$ is norm preserving.\\

If we define the transformation $\Psi: [\ell_{\infty}(F)]_{\Gamma} \rightarrow \ell_{\infty}(F)$ by $u \mapsto y =\Psi u = \sum_{k=1}^{n}\left|(1/k)u_{(t_{1},t_{2})}^{k}\right|$, then we can say that $\Gamma$ is norm preserving.
\end{proof}

\begin{thm}\label{normed}
Integrated and differentiated sequence spaces derived by triangular fuzzy numbers are complete normed sequence spaces with
the norms defined by
\begin{eqnarray}\label{norm}
\|\widehat{u}\|_{\Omega}=\sup_{k}\max\bigg\{\Omega\left|u_{(t_{1},t_{2})}^{k}-v_{(t_{1},t_{2})}^{k}-t_{1}\right|, \Omega\left|u_{(t_{1},t_{2})}^{k}-v_{(t_{1},t_{2})}^{k}\right|,
\end{eqnarray}
\begin{eqnarray*}
\Omega\left|u_{(t_{1},t_{2})}^{k}-v_{(t_{1},t_{2})}^{k}+t_{2}\right|\bigg\}
\end{eqnarray*}
\begin{eqnarray}\label{norm1}
\|\widehat{u}\|_{\Gamma}=\sup_{k}\max\bigg\{\Gamma\left|u_{(t_{1},t_{2})}^{k}-y_{(t_{1},t_{2})}^{k}-t_{1}\right|, \Gamma\left|u_{(t_{1},t_{2})}^{k}-y_{(t_{1},t_{2})}^{k}\right|,
\end{eqnarray}
\begin{eqnarray*}
\Gamma\left|u_{(t_{1},t_{2})}^{k}-y_{(t_{1},t_{2})}^{k}+t_{2}\right|\bigg\},
\end{eqnarray*}
respectively.
\end{thm}

\begin{proof}:
Integrated and differentiated spaces derived by triangular fuzzy numbers are norm isomorphic to
the classical sets consisting of the bounded, convergent and null sequences of fuzzy numbers.
Further, the matrices $\Omega$ and $\Gamma$ are regular. Hence, integrated and differentiated spaces derived by triangular fuzzy numbers are complete
normed spaces with the norms (\ref{norm}), (\ref{norm1}).
\end{proof}

Now, we give the real duals of the integrated and differentiated spaces with triangular fuzzy numbers.\\

We take two spaces of triangular fuzzy valued sequences $X(F)$ and $Y(F)$. We define the set
\begin{eqnarray}\label{dual}
&&T(X(F), Y(F))=\left\{a=(a^{k})\in \omega(F): \left(a^{k}x_{(t_{1},t_{2})}^{k}\right)\in Y(F), \forall x \in X(F)\right\}.
\end{eqnarray}
We denotes the real duals of the sequence spaces $X(F)$ with $\alpha(r)$, $\beta(r)$, $\gamma(r)$. Using the notation (\ref{dual}),
we can define the duals as follows:
\begin{eqnarray*}
&&[X(F)]^{\alpha(r)}=T\big(X(F), \ell_{1}(F)\big),\\
&&[X(F)]^{\beta(r)}=T\big(X(F), cs(F)\big),\\
&&[X(F)]^{\gamma(r)}=T\big(X(F), bs(F)\big),
\end{eqnarray*}
where $cs(F)$ and $bs(F)$ denotes the spaces of convergent and bounded series of fuzzy numbers, respectively.

\begin{lem}\label{lemMTRX}\cite{TaBa1}
The following statements are hold:\\
i. $A\in (\ell_{\infty}(F):\ell_{\infty}(F))$, $A\in (c(F):\ell_{\infty}(F))$, $A\in (c_{0}(F):\ell_{\infty}(F))$
if and only if the condition $\sup_{n}\sum_{k}\overline{d}\left(a_{nk}, \theta\right)<\infty$ holds.\\
ii. $A\in (\ell_{\infty}(F):c_{0}(F))$ if and only if the condition $\lim_{n}\sum_{k}\overline{d}\left(a_{nk}, \theta\right)=0$ holds.\\
iii. $A\in (c_{0}(F):c(F))$ if and only if the conditions $\sup_{n}\sum_{k}\overline{d}\left(a_{nk}, \theta\right)<\infty$ and $\lim_{n}\overline{d}\left(a_{nk}, a^{k}\right)=0 \quad \textit{where} \quad \left(a_{(0,0)}^{k}\right)\in \omega(F)$ hold.\\
iv. $A\in (c_{0}(F):c_{0}(F))$ if and only if the conditions $\sup_{n}\sum_{k}\overline{d}\left(a_{nk}, \theta\right)<\infty$ and $\lim_{n}\overline{d}\left(a_{nk}, a^{k}\right)=0 \quad \textit{where} \quad \left(a_{(0,0)}^{k}\right)\in \omega(F)$ hold with $a^{k}=\theta$ for all $k\in \mathbb{N}$.
\end{lem}

\begin{lem}\label{seng1lemma}\cite{seng1}
Let $K$ is the finite subset of $\mathbb{N}$ and $A$ be an infinite matrix of positive numbers $a_{nk}$.\\
\textbf{i.} $A\in (c(F):c(F))$ if and only if the conditions $\sup_{n}\sum_{k}\overline{d}\left(a_{nk}, \theta\right)<\infty$ and $\lim_{n}\overline{d}\left(a_{nk}, a^{k}\right)=0 \quad \textit{where} \quad \left(a_{(0,0)}^{k}\right)\in \omega(F)\mathbf{}$ hold with $a^{k}=\theta$ for all $k\in \mathbb{N}$.\\
\textbf{ii.} $A\in (c_{0}(F):\ell_{1}(F))$ if and only if $\sup_{K}\sum_{k}\overline{d}\bigg(\sum_{n\in K}a_{n,k},\theta\bigg)<\infty$.
\end{lem}

Let $U=(u_{nk})$ be an infinite matrix and $V=(v_{nk})$ be an inverse matrix of $U$. The matrix $F=(f_{nk})$
defined by $f_{nk}=\sum_{j=k}^{n}a^{j}v_{jk}$, $(0\leq k \leq n)$ and $f_{nk}=0$, $(k>n)$, for all $k,n\in\mathbb{N}$.

\begin{lem}\label{lemBeta}\cite{AltBas}
Let $X$ be a sequence space. Then, \\
\begin{eqnarray*}
&&[X_{U}]^{\beta}=\left\{a=(a^{k})\in \omega: F\in (X:c)\right\}\\
&&[X_{U}]^{\gamma}=\left\{a=(a^{k})\in \omega: F\in (X:\ell_{\infty})\right\}.
\end{eqnarray*}
\end{lem}

\begin{thm}
The sets $d_{1}$ and $d_{2}$ are defined as follows:
\begin{eqnarray*}
&&d_{1}=\left\{a\in \omega(F): \sup_{n\in \mathbb{N}}\sum_{k}\overline{d}\left(\sum_{n\in\mathbb{N}}\sum_{k\in K}\frac{a_{k}}{k}, \theta\right)<\infty\right\}\\
&&d_{2}=\left\{a\in \omega(F): \sup_{n\in \mathbb{N}}\sum_{k}\overline{d}\left(\sum_{n\in\mathbb{N}}\sum_{k\in K}ka_{k}, \theta\right)<\infty\right\}.
\end{eqnarray*}
The $\alpha(r)-$duals of the integrated and differentiated triangular fuzzy sequence spaces are the sets $d_{1}$ and $d_{2}$, respectively.
\end{thm}

The proof of this theorem is obtained from (ii) of Lemma \ref{seng1lemma}.\\

We define the matrices $G=(g_{nk})$ as $g_{nk}=a_{n}/n$, $(1\leq k \leq n)$ and $g_{nk}=0$, $(k>n)$, and $H=(h_{nk})$ as
by $h_{nk}=na_{n}$, $(1\leq k \leq n)$ and $h_{nk}=0$, $(k>n)$ via the sequence $a=a^{i}\in \omega$. Also we define the sets $d_{3}, d_{4}, d_{5}, d_{6}$ as follows:\\
\begin{eqnarray*}
&&d_{3}=\left\{a\in \omega(F): \sup_{n}\sum_{k}\overline{d}(g_{nk}, \theta)<\infty\right\}\\
&&d_{4}=\left\{a\in \omega(F): \lim_{n\rightarrow\infty}\overline{d}(g_{nk}, \theta)  ~ \textrm{exists for each}~ k\in\mathbb{N}\right\}\\
&&d_{5}=\left\{a\in \omega(F): \sup_{n}\sum_{k}\overline{d}(h_{nk}, \theta)<\infty\right\}\\
&&d_{6}=\left\{a\in \omega(F): \lim_{n\rightarrow\infty}\overline{d}(h_{nk}, \theta)  ~ \textrm{exists for each}~ k\in\mathbb{N}\right\}
\end{eqnarray*}

\begin{thm}
$\big\{[X(F)]_{\Omega}\big\}^{\beta(r)}=d_{3}\cap d_{4}$ and $\big\{[X(F)]_{\Gamma}\big\}^{\beta(r)}=d_{5}\cap d_{6}$.
\end{thm}

Using the Lemma \ref{lemMTRX} and Lemma \ref{lemBeta}, the proof of this theorem is obtained. This proof is similar to proof of Theorem \ref{gammathm}.
Then, we only prove the Theorem \ref{gammathm}.\\

\begin{thm}\label{gammathm}
The $\gamma(r)-$duals of the integrated and differentiated triangular fuzzy sequence spaces are the sets $d_{3}$ and $d_{5}$, respectively.
\end{thm}

\begin{proof}
We give only proof of the $\gamma(r)-$dual of the space $[c_{0}(F)]_{\Omega}$. Let $a\in \omega$ and give the matrix $G=(g_{nk})$.
Using the relation (\ref{transfseq1}), we have
\begin{eqnarray}\label{gammadual}
\sum_{k=0}^{n}a^{k}x_{(t_{1},t_{2})}^{k}=\left(Gy_{(t_{1},t_{2})}\right)^{n}
\end{eqnarray}
We obtain that $ax=\big(a^{k}x_{(t_{1},t_{2})}^{k}\big)\in bs(F)$ whenever $x\in [c_{0}(F)]_{\Omega}$ if and only if $Gy_{(t_{1},t_{2})}\in \ell_{\infty}(F)$
whenever $y_{(t_{1},t_{2})}\in c_{0}(F)$ by \ref{gammadual}. From (i) of Lemma \ref{lemMTRX}, we have $\sup_{n}\sum_{k}\overline{d}(e_{nk}, \theta)<\infty$.
Then, we get from the last result that $\big\{[c_{0}(F)]_{\Omega}\big\}^{\gamma}(r)=d_{3}$.\\

Similarly, it is seen that $\big\{[c_{0}(F)]_{\Gamma}\big\}^{\gamma}(r)=d_{5}$ by using the matrix $H=(h_{nk})$.
\end{proof}

\section{Matrix Transformations}

In this section, we characterize the matrix mappings. First, we will give some brevity as follows:\\
\begin{eqnarray}\label{mtrstransf0}
\widetilde{b_{nk}}=\frac{a_{nk}}{k}-\frac{a_{n,k+1}}{k+1} \quad \quad ~ \textrm{or}~ \quad \widetilde{a_{nk}}=\sum_{k=1}^{\infty}\sum_{j=1}^{n}jb_{nj}
\end{eqnarray}
and
\begin{eqnarray}\label{mtrstransf00}
\widehat{b_{nk}}=ka_{nk}-(k+1)a_{n,k+1} \quad \quad ~ \textrm{or}~ \quad \widehat{a_{nk}}=\sum_{k=1}^{\infty}\sum_{j=1}^{n}\frac{1}{j}b_{nj}
\end{eqnarray}

\begin{thm}\label{omegamtrx}
Let $P=(p_{ni})$ and $R=(r_{ni})$ be $\Omega$ matrices and $X(F)$ be any given sequence space.
Then, $P\in \big([c(F)]_{\Omega}:X(F)\big)$ if and only if $(p_{ni})_{i\in\mathbb{N}}\in \ell_{1}(F)$ and $R\in \big(c(F):X(F)\big)$.
\end{thm}

\begin{proof}
We choose the matrices $P=(p_{ni})$ and $R=(r_{ni})$ which are $\Omega$ matrices. Assume that these matrices satisfy the condition (\ref{mtrstransf0}).
Also it is known that the spaces $[c(F)]_{\Omega}$ and $c(F)$ are linearly isomorphic(Theorem \ref{isotheo}).\\

Suppose that $P\in ([c(F)]_{\Omega}:X(F))$. If we take any $z\in c(F)$ then, $R\Omega$ is equal to $P$ and $(p_{ni})_{i\in\mathbb{N}}\in {[c(F)]_{\Omega}}^{\beta(r)}$. Therefore, we obtain that $(r_{ni})_{i\in\mathbb{N}}\in \ell_{1}(F)$ for each $n\in\mathbb{N}$. Hence, $Rz$
exists for each $z\in c(F)$. Then, we have
\begin{eqnarray}\label{mtrstransf1}
\sum_{i}r_{ni}z_{(t_{1},t_{2})}^{i}=\sum_{i}p_{ni}x_{(t_{1},t_{2})}^{i}
\end{eqnarray}
for all $n\in\mathbb{N}$. From the equation \ref{mtrstransf1}, $Rz=Px$ and $R\in (c(F):X(F))$.\\

Conversely, suppose that $R\in (c(F):X(F))$. If we take any $x\in [c(F)]_{\Omega}$, then $Px$ exists. Therefore, we obtain from the following equality as $n\rightarrow \infty$ that $Px=Rz$ and this shows that $P\in ([c(F)]_{\Omega}:X(F))$ and
\begin{eqnarray*}
\sum_{i=0}^{n}p_{ni}x_{(t_{1},t_{2})}^{i}=\sum_{i=0}^{n}\left\{\frac{1}{n}\sum_{i=0}^{n}a_{ni}\right\}z_{(t_{1},t_{2})}^{i}
\end{eqnarray*}
for all $n\in\mathbb{N}$, which is the desired result.
\end{proof}

\begin{thm}\label{gammamtrx}
Let $T=(t_{ni})$ and $Q=(q_{ni})$ be the $\Gamma$ matrices and $X(F)$ be any given sequence space.
Then, $T\in \big([c(F)]_{\Gamma}:X(F)\big)$ if and only if $(a_{ni})_{k\in\mathbb{N}}\in \ell_{1}(F)$ and $Q\in \big(c(F):X(F)\big)$.
\end{thm}

Theorem \ref{gammamtrx} can be proved using (\ref{mtrstransf00}).

\begin{thm}\label{invomegamtrx}
Suppose that the elements of the infinite matrices $\Delta=(\delta_{ni})$ and $\Lambda=(\lambda_{ni})$ are connected with relation
\begin{eqnarray}\label{mtrstransf2}
\lambda_{nk}=\sum_{k=1}^{\infty}\sum_{j=1}^{n}j\delta_{jk}
\end{eqnarray}
for all $n,i\in \mathbb{N}$. Then, $\Delta\in \big(X(F):[c(F)]_{\Omega}\big)$ if and only if $\Lambda\in \big(X(F):c(F)\big)$.
\end{thm}
\begin{proof}
Take $x=\big(x_{(t_{1},t_{2})}^{i}\big)\in X(F)$. Then, we have
\begin{eqnarray*}
\sum_{k=1}^{m}\lambda_{nk}x_{(t_{1},t_{2})}^{i}&=&\sum_{k=1}^{m}\sum_{k=1}^{\infty}\sum_{j=1}^{n}j\delta_{jk}x_{(t_{1},t_{2})}^{i}\\
&=&
\end{eqnarray*}
for all $m,n\in \mathbb{N}$. Then, we obtain $(\Lambda x)_{n}=[\Gamma(\Delta x)]_{n}$ as $m\rightarrow \infty$, for all $n\in \mathbb{N}$.
Therefore, one can observe from here that $\Delta x\in [c(F)]_{\Omega}$ whenever $x\in X(F)$ if and only if $\Lambda x \in c(F)$ whenever $x\in X(F)$.
\end{proof}

\begin{thm}\label{invgmaamtrx}
Suppose that the elements of the infinite matrices $\Pi=(\pi_{ni})$ and $\Sigma=(\sigma_{ni})$ are connected with relation
\begin{eqnarray}\label{mtrstransf3}
\sigma_{nk}=\sum_{k=1}^{\infty}\sum_{j=1}^{n}\frac{1}{j}\pi_{jk}
\end{eqnarray}
for all $n,i\in \mathbb{N}$. Then, $\Pi\in \big(X(F):[c(F)]_{\Gamma}\big)$ if and only if $\Sigma\in \big(X(F):c(F)\big)$.
\end{thm}

\section{Conclusion}
Fuzzy sets have a great progress in every scientific research area. It found many application areas in both theoretical and practical
studies from engineering area to arts and humanities, from computer science to health sciences, and from life sciences to physical
sciences. In the recent years, ordinary fuzzy sets
have been extended to new types and these extensions have been used in many areas such as energy, medicine, material, economics
and pharmacology sciences.\\

The theory of fuzzy sets now encompasses a well organized corpus of basic notions including (and not restricted to) aggregation operations, a generalized theory of relations, specific measures of information content, a calculus of fuzzy numbers.\\

In this paper, using the matrices $\Omega$ and $\Gamma$, we define the new sequence spaces with triangular fuzzy numbers. We compute the real-duals of these spaces
and characterize the matrix classes of these spaces with well-known sequence spaces.\\

The properties and results related to the integrated spaces with triangular fuzzy numbers and differentiated spaces with triangular fuzzy numbers
are more general and more extensive than the corresponding consequences of the classical sets consisting of the bounded, convergent and null
sequences of fuzzy numbers.\\

\end{document}